\documentclass[a4paper,11pt]{article}
\usepackage[paperwidth=21cm,paperheight=29.7cm,left=2.5cm,right=2.5cm,top=2.5cm,bottom=2.5cm]{geometry}
\usepackage{amsmath}
\usepackage[utf8]{inputenc}
\usepackage[ruled]{algorithm2e}
\usepackage{longtable}
\usepackage{varwidth}
\usepackage[labelfont=bf,tableposition=t]{caption}
\usepackage{array}
\usepackage{multirow}
\usepackage{makecell}
\setcellgapes{1pt}
\makegapedcells
\newcolumntype{C}[1]{>{\centering\arraybackslash }b{#1}}
\newcolumntype{I}[1]{%
   >{\centering\arraybackslash}
   p{#1}
   <{}
}%

\usepackage{empheq}
\usepackage[title]{appendix}
\usepackage{comment}

\newcommand{\tab}[1]{\hspace*{#1}}

\newtheorem{defi}{Definition}{\bfseries}{\itshape}
\newtheorem{prop}{Proposition}{\bfseries}{\itshape}
\newtheorem{prop1}{Property}{\bfseries}{\itshape}

\providecommand{\keywords}[1]{\noindent\textbf{\textit{Keywords---}} #1}

\title{Heuristic and exact fixation-based approaches for the discounted 0-1 knapsack problem}

\makeatletter
\renewcommand\@date{{%
  \vspace{-2\baselineskip}%
  \large\centering
  \begin{tabular}{@{}c@{}}
    Christophe Wilbaut\textsuperscript{1,2} \\
    \normalsize christophe.wilbaut@uphf.fr
  \end{tabular}%
  \quad \quad 
  \begin{tabular}{@{}c@{}}
    Raca Todosijevic\textsuperscript{1,2} \\
    \normalsize raca.todosijevic@uphf.fr
  \end{tabular}
  \\\vspace{\baselineskip}
  \begin{tabular}{@{}c@{}}
    Sa\"{i}d Hanafi\textsuperscript{1,2} \\
    \normalsize said.hanafi@uphf.fr
  \end{tabular}
  \quad \quad 
  \begin{tabular}{@{}c@{}}
    Arnaud Fr{\'e}ville\textsuperscript{1}\footnote{``Professeur émérite'' - Now retired}
  \end{tabular}

  \bigskip

  \textsuperscript{1}Univ. Polytechnique Hauts-de-France, LAMIH, CNRS, UMR 8201, F-59313 Valenciennes, France\par
  \textsuperscript{2}INSA Hauts-de-France, F-59313 Valenciennes, France

  \bigskip

  \date{}
}}
\makeatother

\begin{document}

\maketitle

\begin{abstract}
\noindent In this paper we consider the discounted 0-1 knapsack problem (DKP), which is an extension of the classical knapsack problem where a set of items is decomposed into groups of three items. At most one item can be chosen from each group and the aim is to maximize the total profit of the selected items while respecting the knapsack capacity constraint. The DKP is a relatively recent problem in the literature. In this paper we propose a two-phase approach in which the problem is reduced by applying exact and / or heuristic fixation rules in a first phase that can be viewed as a preprocessing phase. The remaining problem can then be solved by dynamic programming. Experiments performed on available instances in the literature show that the fixation techniques are very useful to solve these instances. Indeed, the preprocessing phase greatly reduces the size of these instances, leading to a significant reduction in the time required for dynamic programming to provide an optimal solution. \newline

\keywords{Knapsack problem - Discounted knapsack problem - 
- Fixation - Dynamic Programming}
\end{abstract}
\section{Introduction}
\label{intro}
Knapsack problems arise in many applications from various areas and several variants were derived from the original knapsack problem (KP). The KP is defined as follows. Given a set $N$ of $n$ items, where each item $j\in N$ has associated profit $c_j$ and weight $a_{j}$, the aim of the KP is to select a subset of $N$ in order to maximize the profit of the selected items without exceeding a given capacity $b$ of the knapsack. The KP is among the optimization problems with the simplest linear integer programming formulation since it can be formulated with only one capacity (or resource) constraint. The binary (0-1) variant of the KP may be defined by using binary variables $x_{j}$, associated to each item $j$, that specify if item $j$ is selected ($x_{j} = 1$) or not ($x_{j} = 0$). Then, the 0-1 KP is given as follows:
%
\begin{empheq}[left=$(\text{0-1\,KP})$\empheqlbrace]{align}
    \text{max} & \displaystyle\sum\limits_{j=1}^{n} c_{j}x_{j} \label{KP:obj}\\
    \text{subject to:} & \displaystyle\sum\limits_{j=1}^{n} a_{j}x_{j} \leq b \label{KP:kp}\\
    & x_{j} \in \{0,1\} \tab{2cm} j \in N=\{1,\ldots,n\} \label{KP:bin}
\end{empheq}
%
The KP has been intensively studied from the middle of the twentieth century, starting especially with works of Lorie and Savage \cite{lorie1955}, Gilmore and Gomory \cite{gilmore1966} or Nemhauser and Ulmann \cite{nemhauser1969}. The success of the KP is due to the fact that it can be extended easily to several interesting and challenging variants by adding side constraints (e.g., the multidimensional knapsack problem), by changing the objective (e.g., the quadratic knapsack problem) or by partitioning the set of variables (e.g., knapsack problem with setup) for instance. In fact the KP appears frequently as a sub-problem of other hard optimization problems and a wide range of practical applications can be listed for these family of problems such as cargo loading, cutting stock, capital budgeting or project selection for instance. Several review papers and books were dedicated to the KP and some of its variants (see for instance \cite{salkin1975,kellerer2004,wilbaut2007} among others). 
The discounted 0–1 knapsack problem (DKP) is one of the extensions of the 0–1 KP in which variables are partitioned. It was introduced in 2007 by Guldan in his Master thesis \cite{guldan2007}. The word \textsl{discounted} comes from the promotional discounts activities of real merchants: in the same way as buying two items at the same time during the promotion season, a discount is defined when buying both. In the DKP items are grouped by three and the value of the discounted item is the sum of the other two items in the group. In addition the discounted item consumes less resources than the total depletion of the other two items. Finally, it is allowed to select at most one item  from the group. The aim of the DKP is to maximize the total value associated with the selected items without exceeding the knapsack capacity. 
Formally, the DKP involves a set $M$ of $m$ groups (or classes) of three items denoted by $N_i=\{3i, 3i+1, 3i+2\}$ for $i \in M=\{0, \ldots, m-1\}$, which thus contains items to be possibly packed into the knapsack of capacity $b$. Each item $3i+k$ in group $N_i$, $k \in \{0,1,2\}$ and $i\in M$, is characterized by a profit $c_{3i+k}$ and a weight $a_{3i+k}$ such that
$c_{3i+2} = c_{3i} + c_{3i+1}$ and $\max\{a_{3i}, a_{3i+1}\} < a_{3i+2} < a_{3i} + a_{3i+1}$. The aim of the DKP is to choose at most one item from each group $N_i$ such that the total profit of chosen elements is maximized without having the weight sum to exceed $b$. In the following we use notation $N=\{0,\ldots, 3m-1\}$ to refer to the index set of items whereas $n=|N|=3m$ is the number of items. As in the KP all the coefficients $c_{j}$, $a_j$, $\forall j \in N$ and the capacity $b$ are integer and are assumed to be non negative. 
A linear formulation of the DKP \cite{rong2012} is obtained by introducing three binary variables $x_{3i}, x_{3i+1}$ and $x_{3i+2}$ associated with each group $i \in M$, such that $x_{3i+k}=1$ if and only if the $k^{th}$ item of group $i$ is selected, for $k \in \{0,1,2\}$. Hence, the standard 0-1 linear programming formulation of DKP is given  as follows: \\ 
%
\begin{empheq}[left=$(\text{DKP})$\empheqlbrace]{align}
    \text{max} & \displaystyle\sum\limits_{i \in M} c_{3i}x_{3i} + c_{3i+1} x_{3i+1} + c_{3i+2} x_{3i+2} \label{L1:obj}\\
    \text{s.t.:} & \displaystyle\sum\limits_{i\in M} a_{3i}x_{3i} + a_{3i+1} x_{3i+1} + a_{3i+2} x_{3i+2} \leq b \label{L1:kp}\\
    & x_{3i} + x_{3i+1} + x_{3i+2} \leq 1 \tab{2.6cm} i \in M \label{L1:choice}\\
    & x_{3i}, x_{3i+1}, x_{3i+2} \in \{0,1\} \tab{2.45cm} i \in M \label{L1:bin}
\end{empheq}
%
\par In this model objective function (\ref{L1:obj})  aims to maximize the total  profit of selected items. Constraint (\ref{L1:kp}) is the knapsack constraint related to the $n = 3m$ items in the problem, whereas constraints (\ref{L1:choice}) force the choice of at most one item from each group $i \in M$. The all variables are binary variables as stated in (\ref{L1:bin}).
As mentioned in Rong et al. \cite{rong2012} additional conditions can be imposed to avoid trivial solutions and to guarantee that each item in a given group has a chance to be selected. In particular it is generally assumed that $a_{3i} < a_{3i+1} < a_{3i+2}$, $c_{3i} < c_{3i+1} < c_{3i+2}$, and $a_{3i} + a_{3i+1} > a_{3i+2}$. Finally, it is also assumed that $\displaystyle\sum_{i \in M} a_{3i+2} > b$ and $a_{3i+2} \leq b, i \in M$. 
It may be observed that by deleting constraints (\ref{L1:choice}) the DKP reduces to the KP. As the KP is an NP-hard problem, this implies the DKP is NP-hard as well. 
%
The DKP can be viewed as a particular case of the multiple-choice KP (MCKP) where items are also grouped. However, unlike to the DKP, in the MCKP 
groups do not necessarily contain the same number of  items. In addition specific relations between items in a group are defined in the DKP, which is not the case in the MCKP. \par

The literature on the DKP started with \cite{guldan2007} where Guldan presented an exact algorithm  based on dynamic programming \cite{bellman1957}. Then, Rong et al. \cite{rong2012} presented a natural formulation of the DKP and solved it by combining the core concept used to solve the KP with the dynamic programming. They proposed three variants based on partitioning the problem according to the type of instance and using or not dominance rules. The results showed that the detection of dominated states leads to an increase of the average running time needed to solve an instance. He et al. also used dynamic programming to solve the DKP in \cite{he2016}. Authors derived a recursive formula based on the principle of
minimizing the total weight rather than maximizing the total profit. This new variant of dynamic programming has clearly a lower complexity than the previous one when the sum of profit coefficients is less than the knapsack capacity. They also proposed a fully polynomial-time approximation scheme, a 2-approximation algorithm and a particle swarm optimization (PSO) heuristic for the DKP. The fully polynomial-time approximation scheme is based on the construction of a new instance by modifying only the set of original profit coefficients. This new instance is then solved by a dynamic programming algorithm. Then, authors proposed a greedy algorithm and show that it is a 2-approximation algorithm for the DKP. They demonstrated that the relationship between the ratios of the three items in a given group can be induced by only four cases and they exploited this property in the greedy algorithm. Finally, authors proposed a greedy repair algorithm to deal with infeasible solutions and incorporate it in a PSO algorithm. Some other references related to evolutionary algorithms dedicated to the DKP can be found in the literature. He et al. proposed in \cite{he2016a} two elitist genetic algorithms and two types of greedy strategies to repair and to optimize the individuals. Feng et al. proposed in \cite{feng2018} a multi-strategy monarch butterfly optimization heuristic for solving the DKP. They introduced two effective strategies based on a neighborhood mutation with crowding and a Gaussian perturbation to improve the behaviour of the approach. Then, Feng and Wang designed ten types of moth search in \cite{feng2018a}. Zhu et al. considered in \cite{zhu2017} three differential evolution algorithms for the DKP that mainly differ
in the encoding mechanisms used to represent (feasible) solutions of the problem. The first variant is an adaptation of the hybrid-encoding binary differential evolution algorithm proposed in \cite{he2007} to solve the KP and the SAT problem. Authors also introduced the repair operator proposed in \cite{he2016a} to manage infeasible solutions. This variant mimics the classical formulation
of the DKP with $3m$ binary variables by using an encoding conversion function to transform a $3m$-dimensional real vector into a $3m$-dimensional binary vector. The other two variants exploit another encoding in which a feasible solution is an integer vector in $\{0, 1, 2, 3\}^{m}$ allowing a direct application of the standard differential evolution mechanism. Authors considered two different encoding conversion functions corresponding to the two variants. More recently, He et al. proposed in \cite{he2019} a new kind of method to design an evolutionary algorithm by using algebraic theory. The proposition is mainly based on two evolution operators using the addition, multiplication and inverse operation of the direct product of rings. The first one is called global exploration operator and build a new individual by learning from four different elements that are randomly selected from 
the whole search space. The second one is called local development operator and is based on random changes in a given individual. In the new evolutionary algorithm authors propose to apply these two operators successively and to replace an existing individual in the population if the new one produced by these operations has a better fitness. The approach is adapted to the DKP by using the repair operator proposed in \cite{he2016a}. Very recently, Wu et al. proposed in \cite{wu2020} a discrete hybrid teaching-learning-based optimization algorithm. They introduced self-learning factors in the teacher and learner phases, which balances the exploitation and exploration of the algorithm. They also used to types of crossover to improve the search. The proposed method is evaluated and compared with 7 other population-based algorithms (\cite{he2016,he2016a,feng2018,feng2018a,zhu2017,wu2017}) on a set of 80 available instances in the literature. The results showed that their approach obtained in average the best solutions. \par

In this paper we consider reduction techniques to solve the DKP efficiently. One of the objective of this work is to show that such techniques can be used in a heuristic or an exact process to obtain near optimal or optimal solutions of the input problem in a very short time. More precisely, the main contributions of this work are: ($i$) we propose heuristic and exact fixation rules for the DKP that can be used in a preprocessing phase in order to reduce the size of a given instance; ($ii$) we show that using these techniques is an efficient way to solve all the existing instances in the literature in just a few seconds when applying dynamic programming.

%
The paper is organized as follows. In Section \ref{sec::fix} we present fixation rules that can be applied for solving the DKP. We consider in particular a heuristic fixation process based on the notion of LP-dominance that can be used to accelerate the search and to improve the behavior of heuristic approaches, and an exact fixation procedure based on a group reduction rule which helps to discard an important number of variables in the problem by choosing optimally one item among the three in several groups. Then, we present in Section \ref{sec::exp} the results obtained with our approaches considering the instances available in the literature. Finally, Section \ref{sec::concl} presents some conclusion and perspectives.

\section{Fixation rules}
\label{sec::fix}
Fixation techniques and reduction rules are often used to solve optimization problems in general and knapsack problems in particular. In many cases the aim is to set values for a subset of variables before solving the reduced problem with a dedicated exact method (see e.g., \cite{rong2013}). In other cases a decomposition method is applied for solving more efficiently the initial problem \cite{chen2014,dahmani2016}. Another solution is to solve the reduced problem obtained when fixing a subset of variables with a metaheuristic \cite{wilbaut2006}. Considering the KP several approaches were proposed to reduce the size before applying a branch-and-bound method or a dynamic programming algorithm (see, e.g., \cite{martello1988,pisinger1997}), whereas other implicit enumeration algorithms include a reduction phase (see, e.g., \cite{fayard1975,fayard1982}). In this section we consider the use of both heuristic and exact fixation rules to solve the DKP. In particular, we consider techniques originally proposed for solving the MCKP or other knapsack variants to the DKP case. We start this section with the notion of LP-dominance that can be used to set variables' values in a heuristic way. Then we present a simple reduction process to remove definitively some groups from the initial problem without eliminating any optimal solution of the problem.

\subsection{LP-dominance}
The first fixation rule is adapted from an existing result for the MCKP. It is based on the notion of \textsl{LP-dominance} which can be used to solve efficiently the LP-relaxation of the MCKP \cite{sinha1979}. As mentioned previously the DKP can be viewed as a special case of the MCKP. Then, a given instance $P$ of DKP can be formulated as a special instance $P'$ of MCKP by introducing for each group a dummy item with a null profit and a null weight. More formally, we construct $P'$ from the original DKP instance by considering $m$ groups of four items $N'_i=\{4i, 4i+1, 4i+2, 4i+3\}$ for $i \in M$ such that $c'_{i,k} = c_{3i+k-1}$ and $a'_{i,k} = a_{3i+k-1}$ for $k \in \{1,2,3\}$ and $c'_{i,0} = a'_{i,0} = 0$ $\forall i\in M$. Then, by introducing a binary variable $x_{i,k}$ which takes value 1 if and only if item $k$ in group $N'_i$ (i.e., item $4i+k$) is chosen, for $k \in \{0,1,2,3\}$, the problem instance $P'$ can be formulated as:
\begin{empheq}[left=(P')\empheqlbrace]{align}
    \text{max} & \displaystyle\sum\limits_{i\in M} \sum\limits_{k=0}^{3} c'_{i,k}x_{i,k} \label{MIP2:obj}\\
    \text{s.t.:} & \displaystyle\sum\limits_{i\in M} \sum\limits_{k=0}^{3} a'_{i,k}x_{i,k} \leq b \label{MIP2:kp}\\
    & \sum\limits_{k=0}^{3} x_{i,k} = 1 \tab{1.5cm} i\in M \label{MIP2:choice}\\
    & x_{i,k} \in \{0,1\} \tab{1.4cm} i\in M, k\in \{0,1,2,3\} \label{MIP2:bin}
\end{empheq}
Problem instance $P'$ can be viewed as a reformulation of the previous standard model $\text{DKP}$ where the $m$ groups are explicitly mentioned in the objective function and in the knapsack constraint, leading to a double sum in (\ref{MIP2:obj}) and in (\ref{MIP2:kp}). The particularity of this model comes from the use of the dummy items allowing us to force the choice of exactly one item  from each group (\ref{MIP2:choice}). Please note that in the rest of the paper we use both notations $4i+k$ or $(i,k)$ to refer to item $k$ in group $N'_i$. \par

When dealing with the MCKP the LP-dominance may be used to set some variables at value 0 in an optimal solution of the LP-relaxation. One of the interest of this property is that it can be used as a pre-processing phase since it is only based on the data \cite{sinha1979}. The following definition introduces the LP-dominance when considering instance $P'$ from the original definition for the MCKP.
\begin{defi}
Let P' be the reformulated instance of the DKP and let $i \in M$ be a given group of items. If some items $j,k,l \in N'_{i}$ with $a'_{i,j} < a'_{i,k} < a'_{i,l}$ and $c'_{i,j} < c'_{i,k} < c'_{i,l}$ satisfy 
\begin{eqnarray}
\displaystyle\frac{c'_{i,l}-c'_{i,k}}{a'_{i,l}-a'_{i,k}} \geq \frac{c'_{i,k}-c'_{i,j}}{a'_{i,k}-a'_{i,j}} \label{eq::fix}
\end{eqnarray}
then item $k$ is said to be LP-dominated by items $j$ and $l$.
\end{defi}
From this initial definition we can observe that when dealing with the DKP we know by hypothesis that the profits and the weights of items are always in non-decreasing order into each group, from the fictive item to the last one (i.e., $c'_{i,0} < c'_{i,1} < c'_{i,2} < c'_{i,3}$ and $a'_{i,0} < a'_{i,1} < a'_{i,2} < a'_{i,3}$ for every $i \in M$). In addition only items $(i,1)$ and $(i,2)$ can be dominated in a group. Those observations lead to the following restricted definition.
\begin{defi}
Let P' be the reformulated instance of the DKP and let $i \in M$ be a given group of items. Item (i,1) is LP-dominated if and only if (\ref{eq::fix1_1}) or (\ref{eq::fix1_2}) is satisfied.
\begin{eqnarray}
\displaystyle\frac{c^{'}_{i,2}-c^{'}_{i,1}}{a^{'}_{i,2}-a^{'}_{i,1}} \geq \frac{c^{'}_{i,1}}{a^{'}_{i,1}} \label{eq::fix1_1} \\
\displaystyle\frac{c^{'}_{i,3}-c^{'}_{i,1}}{a^{'}_{i,3}-a^{'}_{i,1}} \geq \frac{c^{'}_{i,1}}{a^{'}_{i,1}} \label{eq::fix1_2} 
\end{eqnarray}
In the same way item (i,2) is LP-dominated if and only if (\ref{eq::fix2_1}) or (\ref{eq::fix2_2}) is satisfied.
\begin{eqnarray}
\displaystyle\frac{c^{'}_{i,3}-c^{'}_{i,2}}{a^{'}_{i,3}-a^{'}_{i,2}} & \geq & \frac{c^{'}_{i,2}-c^{'}_{i,1}}{a^{'}_{i,2}-a^{'}_{i,1}} \label{eq::fix2_1} \\
\displaystyle\frac{c^{'}_{i,3}-c^{'}_{i,2}}{a^{'}_{i,3}-a^{'}_{i,2}} & \geq & \frac{c^{'}_{i,2}}{a^{'}_{i,2}} \label{eq::fix2_2} 
\end{eqnarray}
\end{defi}

From this definition we can apply the following proposition (a proof can be found for instance in \cite{sinha1979}).
\begin{prop}\label{prop::fix}
Let $i \in M$ be a group of items. If item (i,1) (resp. (i,2)) is LP-dominated then an optimal solution to the LP relaxation of P' with $x_{i,1} = 0$ (resp. $x_{i,2} = 0$) exists.
\end{prop}
The characteristics of the DKP allow us not to have to sort the items into the groups, leading to a very short pre-processing phase. Please note that another dominance rule exists to set some variables definitively at their optimal value in an optimal solution of the MCKP (see \cite{sinha1979}). It can be applied if there exist items $(i,j)$ and $(i,j')$ in a group $i$ satisfying $c'_{i,j} \geq c'_{i,j'}$ and $a'_{i,j} \leq a'_{i,j'}$. However, in the case of the DCKP the corresponding propriety cannot be satisfied (here again due to the hypothesis on the profits and the weights).
In this paper we apply Proposition \ref{prop::fix} to set some variables at value 0 in a heuristic way. To achieve this we adapt the algorithm proposed by Zemel \cite{zemel1980} to solve the LP-relaxation of the MCKP. This approach first reformulates the original LP-relaxation of the MCKP into a corresponding LP-relaxation of a KP. Then, this relaxation can be solved efficiently based on the well-known greedy algorithm. Proposition \ref{prop::fix} can be applied while  solving the LP-relaxation. Our method is described in Algorithm \ref{algo::lp}. It consists in the same two main steps. The first one aims at building the KP instance composed by at most $3m$ items (since we can eliminate some variables according to Proposition \ref{prop::fix}). In the KP instance we associate with each item $k \in \{1,2,3\}$ in group $i \in M$ of the DKP a profit noted $c''_{i,k}$ and a weight $a''_{i,k}$. The procedure to build the KP was initially presented and justified by Zemel in \cite{zemel1980}. In the case of the DKP it can be summarized as follows. Based on the fact that items are already sorted in each group according to their weight in the knapsack constraint, the fictive item is first eliminated in all the groups. Then, when considering a group composed only by items that are not LP-dominated, equations (\ref{eq::trans1}-\ref{eq::trans2}) and (\ref{eq::trans3}-\ref{eq::trans4}) are applied respectively to compute the profit and weight values in the KP instance:
\begin{align}
c^{''}_{i,1} & = c^{'}_{i,1} \label{eq::trans1} \\
c^{''}_{i,k} & = c^{'}_{i,k} - c^{'}_{i,k-1} \hspace*{0.5cm} k = 2, 3 \label{eq::trans2} \\
a^{''}_{i,1} & = a^{'}_{i,1} \label{eq::trans3} \\
a^{''}_{i,k} & = a^{'}_{i,k} - a^{'}_{i,k-1} \hspace*{0.5cm} k = 2, 3 \label{eq::trans4}
\end{align}
The incremental profit $c^{''}_{i,k}$ in group $k$ is a measure of how much we gain if item $k$ is chosen instead of item $k-1$. The incremental weight $a^{''}_{i,k}$ has a similar interpretation. One can observe that if a given item is LP-dominated then the procedure is adapted to discard this value as shown in Algorithm \ref{algo::lp}.
In this algorithm we use the following convention: an LP-dominated item has its profit and weight fixed to 0 in the KP instance.
The step consisting of eliminating the LP-dominated items and building the KP instance is described between line \ref{algo::lp:start_step1} and line \ref{algo::lp:end_step1}. The LP-dominated items are added into a set denoted to as $F^{0}$ (and defined in line \ref{algo::lp:init}) which thus contains variables that can be set at 0 in an optimal solution of the LP-relaxation.

In the second step of the algorithm (from line \ref{algo::lp:start_step2} to line \ref{algo::lp:end_step2}) the LP-relaxation of the KP instance is solved with a greedy algorithm when items are ordered according to their incremental efficiencies denoted to as $e''_{i,k}$ in Algorithm \ref{algo::lp}. By convention the efficiency of an LP-dominated item is set to $-\infty$ so that it would not be added in the solution (see lines \ref{algo::lp:eff1} and \ref{algo::lp:eff2}). During this solving a feasible solution is also obtained by discarding the fractional item in the LP solution. Line 19 is used to avoid the selection of more than one item in a given group, whereas lines \ref{algo::lp:2frac} to \ref{algo::lp:end_step2} manage the case where the LP solution has two fractional items. In that case they are necessarily in the same group (see Proposition 3 in \cite{zemel1980}). 
The last part of the algorithm (from line \ref{algo::lp:start_step3} to line \ref{algo::lp:end_step3}) can be used to try adding a few items in the current feasible solution $x$ in a greedy way, where we only consider non LP-dominated items and groups where no item were previously selected.
Variables returned by Algorithm \ref{algo::lp} in set $F^{0}$ can be set definitively in the problem if we want to use the fixation as a heuristic.

\resizebox{0.8\textwidth}{!}{%
\begin{algorithm}[H]
\SetNoFillComment
\caption{Solving the LP-relaxation and building a greedy solution}\label{algo::lp}
\hspace*{10pt}  \parbox{12in} {{\bf Function} {\tt LP-Greedy} \par
 \nl $F^{0} \leftarrow \oslash$ \; \label{algo::lp:init}
 \tcc{Step1: Eliminate LP-dominated items and build an instance of KP}
 \nl \For {$i$ from 1 to $m$}{ \label{algo::lp:start_step1}
    \tcc{Check if item $(i,1)$ is LP-dominated}
 \nl    \eIf{equation (\ref{eq::fix1_1}) or equation           (\ref{eq::fix1_2}) is satisfied}{
    \nl     $c''_{i,1} = a''_{i,1} \leftarrow 0$ ; $e''_{i,1} \leftarrow -\infty$; \label{algo::lp:eff1}
        $F^{0} \leftarrow F^{0} \cup {(i,1)}$ \;
    }
    {
    \nl     $c''_{i,1} \leftarrow c'_{i,1}$ ; $a''_{i,1} \leftarrow a'_{i,1}$ ; $e''_{i,1} \leftarrow c''_{i,1} / a''_{i,1}$ ;
    }
    \tcc{\footnotesize Check if item $(i,2)$ is LP-dominated}
    \nl \eIf{equation (\ref{eq::fix2_1}) or equation           (\ref{eq::fix2_2}) is satisfied}{
    \nl     $c''_{i,2} = a''_{i,2} \leftarrow 0$ ; $e''_{i,2} \leftarrow -\infty$; \label{algo::lp:eff2}
        $F^{0} \leftarrow F^{0} \cup {(i,2)}$ \;
    \nl $c''_{i,3} \leftarrow c'_{i,3} - c''_{i,1}$ ; $a''_{i,3} \leftarrow a'_{i,3} - a''_{i,1}$ ; $e''_{i,3} \leftarrow c''_{i,3} / a''_{i,3}$ ;
    }
    {
    \nl     $c''_{i,2} \leftarrow c'_{i,2} - c''_{i,1}$ ; $a''_{i,2} \leftarrow a'_{i,2} - a''_{i,1}$ ; $e''_{i,2} \leftarrow c''_{i,2} / a''_{i,2}$ \;
    \nl $c''_{i,3} \leftarrow c'_{i,3} - c''_{i,2}$ ; $a''_{i,3} \leftarrow a'_{i,3} - a''_{i,2}$ ; $e''_{i,3} \leftarrow c''_{i,3} / a''_{i,3}$ ; \label{algo::lp:end_step1}
    }
 }
 \tcc{Step 2: Solve the KP instance in a greedy way}
 \nl $\overline{v} = \underline{v} \leftarrow 0$ ; \label{algo::lp:start_step2}
 $\overline{b} \leftarrow b$; $j \leftarrow 1$ \;
 \nl $\bar{x}_{i,k} = x_{i,k} \leftarrow 0$, $\forall i\in M, k\in \{1,2,3\}$\;
 \nl Sort the items according to non-increasing order of $e''_{i,k}$ values, $\forall i\in M, k \in \{1,2,3\}$\; 
 \nl \While{$\overline{b} > 0$}{
    \nl Let $(i',k')$ the original indexes of next item in the order \;
    \nl \eIf{$\overline{b} > a''_{i',k'}$}{
        \nl $\overline{v} \leftarrow \overline{v} + c''_{i',k'}$ ; $\overline{b} \leftarrow \bar{b} - a''_{i',k'}$\;
        \nl $\bar{x}_{i',k'} = x_{i',k'} \leftarrow 1$\;
        \nl $\bar{x}_{i',k} = x_{i',k} \leftarrow 0, \forall k \in \{1,2,3\}$ such that $k \neq k'$ \;
    }
    {   
        \nl $\underline{v} \leftarrow \overline{v}$ \; \label{algo::lp:2frac}
        \nl $\bar{x}_{i',k'} \leftarrow \overline{b} / a''_{i',k'}$ ; $\overline{v} \leftarrow \overline{v} + c''_{i',k'}\bar{x}_{i',k'}$ ; $\bar{b} \leftarrow 0$ \;
        \nl \If{$\exists k \neq k'$ in group $i'$ such that $\bar{x}_{i',k} = 1$}{
            \nl $\bar{x}_{i',k} \leftarrow 1 - \bar{x}_{i',k'} $; 
            $x_{i',k} \leftarrow 0$ \;
        }
    }
    \nl $j \leftarrow j + 1$ \; \label{algo::lp:end_step2}
 }
 \tcc{Step 3: Try to fill the feasible solution}
 \nl \While{$j \leq 3m - \left|F^{0}\right|$}{ \label{algo::lp:start_step3}
    \nl  Let $(i',k')$ the original indexes of next item in the order \;
    \nl \If{$\overline{b} > a''_{i',k'}$ and $x_{i',k} = 0, \forall k \in \{1,2,3\}$}{
        \nl $x_{i',k'} \leftarrow 1$ ;
        $\underline{v} \leftarrow \underline{v} + c''_{i',k'}$ ; $\overline{b} \leftarrow \overline{b} - a''_{i',k'}$ \;
    }
    \nl $j \leftarrow j + 1$ \; \label{algo::lp:end_step3}
 }
 \nl Return $(\bar{x}, \overline{v})$, $(x, \underline{v})$ and $F^{0}$\;
}
\end{algorithm}
}

As we will show in Section \ref{sec::exp} devoted to the computational experiments this heuristic can be used for almost all the instances of the DKP from the literature without discarding an optimal solution, leading to very fast near optimal solutions. In addition to set $F^{0}$ this algorithm provides an optimal solution $\overline{x}$ of the LP-relaxation of the DKP, the associated upper bound $\overline{v}$, a feasible solution $x$ and the corresponding lower bound $\underline{v}$.

\subsection{Reducing the number of groups}

In the previous section we presented a technique that can be used to fix variables in a heuristic way. In this section we consider another technique to set variables definitively at their optimal value in an optimal solution of the DKP. Thus in that case the fixation is valid for the original problem. 
This fixation rule comes from a well-known property often used when solving variants of knapsack problems. Let $P$ be an instance of the DKP and let $(P\left|x_{i,k} = \alpha\right.)$ be instance $P$ in which we fix only variable $k$ in group $i$ at $\alpha \in \{0, 1\}$, for $i \in M$ and $k \in \{1,2,3\}$. In addition, let $\overline{v}(P\left|x_{i,k}=\alpha\right.)$ be an upper bound of problem $(P\left|x_{i,k}=\alpha\right.)$ and $\underline{v}$ be a lower bound on the optimal value of the original problem $P$. Finally, let $y$ be the feasible solution associated with $\underline{v}$. We can then apply the following Property \ref{prop::pfix}, where notation $\lfloor \gamma \rfloor$ is used to refer to the rounding function, which returns the maximum integer number not greater than $\gamma$.
\begin{prop1}\label{prop::pfix}
If $\lfloor\overline{v}(P\left|x_{i,k}=\alpha\right.)\rfloor \leq \underline{v}$, then there exists an optimal solution $x^{*}$ of problem $P$ where $x^{*}_{i,k} = 1 - \alpha$ or solution $y$ is optimal for problem $P$.
\end{prop1}
Our idea is to exploit the optimal solution $\bar{x}$ of the LP-relaxation of the DKP provided by Algorithm \ref{algo::lp} as a starting point. In fact, we do not consider all the variables in the fixation process. Indeed, when considering the DKP the most interesting result is the fixation of one variable in a group since that allows to remove this group since the other two variables are automatically set to value 0. In addition, according to the special structure of the DKP and the coefficients in the problem it is clear that in a given group $i$ variable $x_{i,3}$ is the one for which the fixation is most probable. Thus, we apply Property \ref{prop::pfix} only with $\alpha = 0$ and for every item $(i,3)$, $i\in M$ such that $\bar{x}_{i,3} = 1$. Indeed, if $\bar{x}_{i,3} = 0$ for a given $i$ then the value $\lfloor\overline{v}(P\left|x_{i,3}=\alpha\right.)\rfloor = \lfloor\bar{v}\rfloor$. Even if in practice most of items $(i,3)$ are chosen in solution $\overline{x}$ in the worst case we thus need to solve $m$ LP-relaxations. \par

Solving the LP-relaxation of problem $(P\left|x_{i,3}=0\right.)$ for a given group $i$ requires some adjustments in steps 1 and 2 of Algorithm \ref{algo::lp}. Indeed, some LP-dominated items in the original problem are not yet dominated when item $(i,3)$ is eliminated (when variable $x_{i,3}$ is set to 0), thus implying some modifications in the construction of the knapsack instance and in the order of items when considering the efficiency measure. However, it is not necessary to restart from scratch and step 1 can be almost avoided since only items $(i,1)$ and $(i,2)$ can be impacted. Once the KP instance is built its LP-relaxation can be solved with the same principle as in Algorithm \ref{algo::lp} to provide both the upper bound $\overline{v}(P\left|x_{i,3}=0\right.)$ and a lower bound for problem $(P\left|x_{i,3}=0\right.)$ which is also valid for the initial DKP. We use notations $\overline{v}_{i}$ (resp. $\underline{v}_{i}$) for this upper (resp. lower) bound and notation $xg_{i}$ for the greedy solution corresponding to this lower bound. Thus, during the process we generate several greedy solutions, in such a way we may improve our initial feasible solution and improve the number of variables that can be set. The procedure is summarized in Algorithm \ref{algo::lps}.

In this algorithm we consider as inputs the results of Algorithm \ref{algo::lp}, in particular the optimal solution $\overline{x}$ of the LP-relaxation of the DKP, the initial lower bound $\underline{v}$, the set $F^{0}$ to know LP-dominated items. To simplify the presentation we also suppose that the KP instance (i.e., vectors $c''$, $a''$ and $e''$) built to solve the LP-relaxation remains available. We use the following notations. Set $F^{1}$ refers to the set of groups that can be eliminated by fixing the corresponding last item at value 1. Thus, we just store in this set the index $i$ of such a group if necessary. Notation $\underline{v}^{best}$ (resp. $x^{best}$) is used to refer to the best lower bound (resp. feasible solution) during the process. 
As explained above the main loop in the algorithm explores the $m$ sub-problems restricted to those where $\overline{x}_{i, 3} = 1$. Between lines \ref{algo::lps::start_adapt} and \ref{algo::lps::end_adapt} we first modify the KP instance by checking if items that were previously LP-dominated are always dominated. In fact in our case when considering group $i$ at the current iteration item $(i,2)$ cannot be LP-dominated since it is the last item in the group. In addition item $(i,1)$ can be only dominated by item $(i,2)$, thus only equation \ref{eq::fix1_1} can be satisfied. 

Then, we apply steps 2 and 3 of Algorithm \ref{algo::lp} to obtain the current upper (resp. lower) bound $\overline{u}_i$ (resp. $\underline{v}_{i}$) and the corresponding feasible solution (lines \ref{algo::lps::step2} and \ref{algo::lps::step3}). Once we have the bounds we then can check if a better lower bound (and feasible solution) has been obtained and if the current group $i$ can be definitively fixed according to Property \ref{prop::pfix} between lines \ref{algo::lps::start_check} and \ref{algo::lps::end_check}. The algorithm returns the set $F^{1}$, the best lower bound and the corresponding best feasible solution. \par

\resizebox{0.8\textwidth}{!}{%
\begin{algorithm}[H]
\SetNoFillComment
\caption{Reducing the number of groups in the DKP}\label{algo::lps}
\hspace*{10pt}  \parbox{12in} {{\bf Function} {\tt UB-Fix} \par
 \nl $F^{1} \leftarrow \oslash$ ; $\underline{v}^{best} \leftarrow \underline{v}$ \;
 \nl \For {$i$ from 1 to $m$}{
    \nl \If{$\overline{x}_{i,3} = 1$}{
    \tcc{Adapt the construction of the KP instance only for group $i$}
        \nl Save the values of $c''_{i,k}, a''_{i,k}$ and $e''_{i,k}$ in the KP instance associated with $\overline{x}$, $\forall k \in \{1,2,3\}$ \; \label{algo::lps::start_adapt}
        \nl \eIf{$(i,1) \in F^{0}$}{
        \nl     \eIf{equation \ref{eq::fix1_1} is satisfied}{
        \nl         $c''_{i,1} = a''_{i,1} \leftarrow 0$ ; $e''_{i,1} \leftarrow -\infty$ \;
            }
            {
        \nl         $c''_{i,1} \leftarrow c'_{i,1}$ ; $a''_{i,1} \leftarrow a'_{i,1}$ ; $e''_{i,1} \leftarrow c''_{i,1} / a''_{i,1}$ ;
            }
        }
        {
        \nl     $c''_{i,1} \leftarrow c'_{i,1}$ ; $a''_{i,1} \leftarrow a'_{i,1}$ ; $e''_{i,1} \leftarrow c''_{i,1} / a''_{i,1}$ \;
        \nl     $c''_{i,2} \leftarrow c'_{i,2} - c''_{i,1}$ ; $a''_{i,2} \leftarrow a'_{i,2} - a''_{i,1}$ ; $e''_{i,2} \leftarrow c''_{i,2} / a''_{i,2}$ \;
        }
        \nl $c''_{i,3} \leftarrow 0$ ; $a''_{i,3} \leftarrow 0$ ; $e''_{i,3} \leftarrow -\infty$ \; \label{algo::lps::end_adapt}
    \tcc{Solve the LP and obtain the upper and lower bounds}
        \nl Apply Steps 2 and 3 of Algorithm \ref{algo::lp} to solve the current KP instance \; \label{algo::lps::step2}
        \nl Obtain $\overline{v}_{i}$, $\underline{v}_{i}$ and $xg_{i}$ \; \label{algo::lps::step3}
        \nl \If{$\underline{v}_{i} \geq \underline{v}^{best}$}{ \label{algo::lps::start_check}
        \nl $\underline{v}^{best} \leftarrow \underline{v}_{i}$ ;
            $x^{best} \leftarrow xg_{i}$ \;
        }
        \nl \If{$\lfloor\overline{v}_{i}\rfloor \leq \underline{v}^{best}$}{
        \nl $F^{1} \leftarrow F^{1} \cup \{i\}$ \; \label{algo::lps::end_check}
    }
    \nl Restore values $c''_{i,k}, a''_{i,k}$ and $e''_{i,k}$ $\forall k \in \{1,2,3\}$ \;
    }
    \nl Return $F^{1}$, $\underline{v}^{best}$, $x^{best}$ \;
 }
}
\end{algorithm}
}

Fixation rules described in this section can be used in different ways. For example it is possible to fix LP-dominated variables according to Proposition \ref{prop::fix} and then to solve the reduced problem. In that case we say that the fixation is heuristic since it is valid when solving the LP-relaxation of the DKP only. Another strategy consists in applying Algorithm \ref{algo::lps} and Property \ref{prop::pfix} to fix definitively some variables $x_{3i+2}$ and thus to eliminate the corresponding groups. In that case we say the fixation is optimal since this property is valid when solving the original DKP.

In our experiments we consider the use of these two techniques with dynamic programming (DP) to solve the instances from the literature.
Dynamic programming was already applied on the DKP in \cite{rong2012} and in \cite{he2016}. Rong et al. presented in \cite{rong2012} a natural extension of the sequential DP algorithm from the classical 0-1 KP to the DKP.
In that case the principle is to look for maximizing the total profit value with the given sum of weight coefficients. Later, He et al. proposed in \cite{he2016} a variant in which the principle is to minimize the total weight with the given sum of value coefficients. The principle of both methods is quite similar and their complexity depends mainly on the size of the problem (in particular the value $b$ or the value $\sum\limits_{i\in M}c_{3i+2}$). 
In this paper we use the natural recursion formula based on the capacity of the knapsack $b$. For the DKP the DP process is based on $m$ stages where each stage corresponds to a set of three variables, where at most one variable can be set to value one. This property allows the use of a similar structure of the DP algorithm to that for solving the knapsack problem. Then, a stage is made up of $b+1$ states defined by the recursive equation according the data of the problem. For the sake of clarity we reintroduce this recursion formula \cite{rong2012}. Let us denote by $v_{i}(\beta)$ the objective value of state $\beta \in \{0,\ldots,b\}$ at stage $i \in \{0,\ldots,n-1\}$. When $i=0$ we have :
\begin{equation}
	{v_0(\beta)=}\left\lbrace
		\begin{aligned}
			0	& \text{ if } 0 \leq \beta < a_0 \\
			c_0	& \text{ if } a_0 \leq \beta < a_1 \\
			c_1 & \text{ if } a_1 \leq \beta < a_2 \\
			c_2 & \text{ if } a_2 \leq \beta \leq b
		\end{aligned}
	\right.
\end{equation}
Then, values $v_{i}(\beta)$ associated with the following stages $1, \ldots, m-1$ and $\beta = 0, \ldots, b$ can be defined as follows:
\begin{equation}
	{v_i(\beta)=}\left\lbrace
		\begin{tabular}{ll}
\hspace*{-7.7cm}$v_{i-1}(\beta)$	& if $0 \leq \beta < a_{3i}$ \\
$\max \{ v_{i-1}(\beta), v_{i-1}(\beta - a_{3i}) + c_{3i} \}$	& if $a_{3i} \leq \beta < a_{3i+1}$ \\
$\max \{ v_{i-1}(\beta), v_{i-1}(\beta - a_{3i}) + c_{3i}, v_{i-1}(\beta - a_{3i+1}) + c_{3i+1} \}$ &  if $a_{3i+1} \leq \beta < a_{3i+2}$ \\
$\max \{ v_{i-1}(\beta), v_{i-1}(\beta - a_{3i}) + c_{3i}, v_{i-1}(\beta - a_{3i+1}) + c_{3i+1},$ & \\
\hspace*{0.8cm}$v_{i-1}(t-a_{3i+2}) + c_{3i+2} \}$ & if  $a_{3i+2} \leq \beta \leq b$
		\end{tabular}
	\right.
\end{equation}
The optimal value corresponds to $v_{m-1}(b)$.

It is important to note that in the rest of the paper when we use the DP to solve a given instance we use a classical implementation of the previous recursive equations and we do not implement any optimization or dominance techniques. For instance we use the traditional way to find an optimal solution by applying a backtracking phase through the set of states. This technique needs $m \times (b+1)$ memory space to save a variable value (i.e., an integer to know which item is selected from the group associated with the stage) for each state corresponding to each stage. Our objective is to show that combining the techniques based on the fixation of variables with DP is an efficient approach to solve the DKP.

\section{Computational results}
\label{sec::exp}
This section is devoted to the presentation of the computational results obtained with our approaches. All the algorithms presented in this paper were coded in {\tt C++} language, implemented and compiled with {\tt Visual Studio} tools on a Windows 10 platform. The tests were carried out on an {\tt HP EliteBook} with 8GB of RAM and an {\tt Intel CORE i7} processor with 2.60GHz. We used {\tt -O2 -Oi -Ot} options in {\tt Visual Studio} to optimize the code and to further accelerate the program. All the CPU times reported in this section were obtained using {\tt clock} function and {\tt CLOCKS\_PER\_SEC} macro.

Rong et al. \cite{rong2012} provided  some characteristics and dominance rules between items according to the correlation of the data and the difficulty of the instance. Then they proposed a generator for three types of instances: uncorrelated, weakly correlated and strongly correlated, respectively. The results obtained by Rong et al. \cite{rong2012} showed that: i) correlated instances are easier to solve than uncorrelated instances for the DKP and ii) strongly correlated instances are easier to solve  than weakly correlated instances. It can be observed that the instances used in most of the papers available in the literature are those provided later by He et al. \cite{he2016}. He et al. proposed a few adjustments to the generator of Rong et al. and developed a fourth type of instances, so-called  inverse strongly correlated instances. All available instances are organized into two data sets (or groups), each one containing 40 instances: 10 instances for each correlation type with $n$ varying between 100 and 1000. These instances are used in this paper to asses performances of the proposed algorithms. 
In the following we refer by the percentage gap, of a solution returned by certain algorithm from the optimal one on a given instance, the value computed as  $100\times\frac{opt-lb}{opt}$, where $opt$ is the optimal value of the problem and $lb$ is the value of the solution. We start this section by solving Linear Programming (LP) relaxation and assessing the ability of fixation rules to reduce the size of the DKP instances in Section \ref{sec::exp_lp}. 
Then Section \ref{sec::exp_dp} demonstrates how fixation rules can also drastically improve the behavior of the dynamic programming for solving the DKP.

\subsection{Solving the LP-relaxation and fixing variables}
\label{sec::exp_lp}
In this section we examine the results obtained when solving the LP-relaxation by Algorithm \ref{algo::lp} and Algorithm \ref{algo::lps}. The aim of these tests are mainly: ($i$) to provide information about the gap between the optimal value of the problem and the LP-bound ; ($ii$) to evaluate the quality of the lower bound derived from the resolution of the LP-relaxation ; ($iii$) to show that the iterative process associated with the solving of the $m$ (at most) upper bounds can be useful to derive a stronger lower bound ; ($iv$) to evaluate capability of both the LP-dominance fixation rule and the reduction process, to decrease size of an instance. A synthesis of the results is presented in Table \ref{tab:lp}, while the detailed results for the 80 instances are provided in 
Tables \ref{tab:detailed_lp_set1} and \ref{tab:detailed_lp_set2} in Appendix \ref{app:lp}. 
In Table \ref{tab:lp} and in the following tables in the paper we report the average results over all instances having the same  correlation type and belonging to the same group: notation {\tt Unc.} (resp. {\tt Weak.}, {\tt Strong.}, {\tt Inv.}) refers to the average results over the 10 uncorrelated (resp. weakly correlated, strongly correlated, inverse strongly correlated) instances of set (or group) 1 or 2. We also provide the average results over all  80 instances. In Table \ref{tab:lp}, columns 2-4 report the results obtained by solving  the LP-relaxation : column {\tt $\lfloor\text{LP}\rfloor$-opt} gives the average difference between the LP value, rounded down, and the optimal value of the problem; column {\tt opt-LPg} corresponds to the average difference between the optimal value and the lower bound value of a greedy solution obtained by simply rounding solution of the LP-relaxation; Column {\tt LP\_Dom} provides the percentage of variables that are LP-dominated according to Proposition \ref{prop::fix}. 
Then, the results of the reduction phase are provided in columns 6 to 8, where {\tt opt-lb} corresponds to the difference between the optimal value and the final lower bound returned by Algorithm \ref{algo::lps}. Columns {\tt Red.} and {\tt Red.+LP\_Dom} refers to the percentage of classes that can be definitively eliminated from the problem since the corresponding $x_{3i+2}$ variable can be set at 1 and the total percentage of variables that can be set at 1 or 0 if we combine the two methods, respectively. Finally, columns {\tt CPU} report the CPU time in seconds needed by the procedures.

\begin{table}[htb]
\caption{Results for the LP-relaxation and reduction process.} \label{tab:lp}
\begin{center}
\begin{scriptsize}
\begin{tabular}{|I{2.2cm}|I{1.2cm}|I{1.2cm}|I{1.2cm}|I{1cm}|I{1.1cm}|I{1.1cm}|I{1.2cm}|I{1cm}|}
\cline{2-9}\multicolumn{1}{c|}{} & \multirow{2}{*}{\textbf{$\lfloor\text{LP}\rfloor$-opt}} & \multirow{2}{*}{\textbf{opt-LPg}} & \multirow{2}{*}{\textbf{LP\_Dom}} & \multirow{2}{*}{\textbf{CPU}} & \multirow{2}{*}{\textbf{opt-lb}} & \multirow{2}{*}{\textbf{Red.}} & \textbf{Red.+} & \multirow{2}{*}{\textbf{CPU}} \\
\multicolumn{1}{c|}{} & & & & & & & \textbf{LP\_Dom} & \\
\hline
Unc.  & 42.80 & 140.40 & 46.90 & 0.004 & 41.50 & 58.72 & 76.07 & 0.039 \\
\hline
Weak. & 35.20 & 280.70 & 64.62 & 0.005 & 28.00 & 59.20 & 84.35 & 0.051 \\
\hline
Strong. & 21.50 & 409.70 & 60.43 & 0.005 & 38.80 & 62.97 & 81.91 & 0.071 \\
\hline
Inv.  & 39.60 & 432.40 & 64.78 & 0.004 & 36.50 & 59.93 & 84.76 & 0.047 \\
\hline
Overall Group 1 & 34.78 & 315.80 & 59.18 & 0.004 & 36.20 & 60.20 & 81.77 & 0.052 \\
\hline
Unc.  & 10.30 & 44.50 & 45.41 & 0.006 & 9.90  & 59.82 & 76.19 & 0.058 \\
\hline
Weak. & 5.20  & 53.00 & 58.42 & 0.007 & 7.00  & 60.20 & 80.18 & 0.053 \\
\hline
Strong. & 6.40  & 266.80 & 52.70 & 0.004 & 16.00 & 61.61 & 77.56 & 0.069 \\
\hline
Inv.  & 8.40  & 40.80 & 59.83 & 0.008 & 4.80  & 53.90 & 77.79 & 0.035 \\
\hline
Overall Group 2 & 7.58  & 101.28 & 54.09 & 0.006 & 9.43  & 58.88 & 77.93 & 0.054 \\
\hline
Overall 80 inst. & 21.18 & 208.54 & 56.64 & 0.005 & 22.81 & 59.54 & 79.85 & 0.053 \\
\hline
\end{tabular}%
\end{scriptsize}
\end{center}
\end{table}

Several interesting observations can be made from the values reported in Table \ref{tab:lp}. First, values in column {\tt $\lfloor\text{LP}\rfloor$-opt} show that the difference between the LP-relaxation and the optimal value of the problem is very small. That is true for all the instances, even more for the second set of instances. 
%
Then, even if the gap between the LP bound and the optimal value of the problem is very tight, the rounding LP-solution does not necessarily yield a solution which is as close to an optimal solution of the problem as it could be expected. However, values reported in column {\tt opt-Lpg} show that extending this rounding solution with a simple greedy approach lead in general to good feasible solutions. To be more precise the gap between this lower bound and the optimal value of the problem is bounded between 0\% (for instance name {\tt wdkp2\_6}) in the best case and 0.49\% in the worst case (instance name {\tt sdkp2\_1}), whereas the average deviation is 0.04\%.
The CPU time needed to solve the LP-relaxation and to obtain a lower bound is negligible, as observed from column {\tt CPU}. Regarding the quality of the lower bound we can observe from column {\tt opt-lb} that the iterative process described in Algorithm \ref{algo::lps} provides in general much better solutions. In particular, this approach leads to 3 and 7 optimal solutions for instances in groups 1 and 2, respectively, although it is not able to improve the initial solution for all the instances. The worst gap is 0.025\% for instance {\tt skdp1\_1}, whereas the average gap falls to 0.004\% which is very tight. A very interesting point is that the CPU time associated with this algorithm is always clearly less than 1 second, even for instances with $n = 1000$. Thus we may conclude that this approach is a very efficient heuristic to build a feasible solution for the DKP in general.
Considering the fixation, column {\tt LP-Dom} indicates that the total number of LP-dominated variables in the LP-relaxation is not insignificant since it is generally superior to 50\%, except for uncorrelated instances. 
Values in column {\tt Red.} are probably more important since they correspond to the percentage of variables that can be definitively eliminated from the initial problem by setting the corresponding $x_{3i+2}$ variables to 1. 
These values are always between 45\% and 72\% with an average value on 59.5\%. Finally, if we combine both the exact and the heuristic fixation rules then we can observe from column {\tt Red.+LP-Dom} that the percentage of variables that can be fixed ranges from 68\% up to 89\%, with an average percentage of 80\%.

\subsection{Combining dynamic programming with fixation}
\label{sec::exp_dp}
In this section we apply dynamic programming after using the reduction techniques.
%
%
%
As mentioned previously in our experiments we used the natural recursion formula based on the capacity of the knapsack $b$. We implemented two versions of DP algorithms when considering or not an optimal solution. Both experiments provide information about the impact of the fixation processes. The results are reported in Table \ref{tab::dp}, while the detailed results over the 80 instances are given in 
Tables \ref{tab::dp_det_set1} and \ref{tab::dp_det_set2} in Appendix \ref{app:dp}.
In Table \ref{tab::dp} we report in columns {\tt Init.Pb.} the CPU times needed to solve the initial problem with our implementation of the DP algorithm. Then, columns {\tt LP-Dom} contain the running times of the same algorithm when it is applied on the reduced problem after applying the LP-dominance rule only (the running time of the fixation process (see Table \ref{tab:lp}) is included in overall running time) and the difference between the optimal value and the value obtained by the DP algorithm. 
Then, column {\tt Red.+LP-Dom} provides the running times when we combine the class reduction phase with the LP-dominance. In that case the final objective value is necessarily the same as in the {\tt LP-Dom} case. Finally, column {\tt Red.} contain the running times when the DP algorithm is applied on the reduced problem obtained when applying only the class reduction. In that case the final value is necessarily the optimal value of the initial problem. \par

Several interesting observations can be made from Table \ref{tab::dp}. First, when considering the initial problem we can see that it is possible to solve instances of set 2 in less than 2 seconds in average, even if we consider the building of the final solution. Detailed experiments showed that for instance set 1, when considering the backtracking phase, values are very similar for medium size instances but logically grow with $n$. In fact, 
the running time needed to solve instances of set 1 is still reasonable for an exact approach since it never exceeds one minute. The difference between the running times needed for solving instances  
in both sets can be explained by the generator used to produce the instances 
(interesting readers are invited to consult \cite{he2016} for more details.) In practice high values for coefficient $b$ can be a limit for the DP algorithm according to the resources available on the machine test. It was for example impossible for us to solve all the instances of set 1 in a {\tt WIN32} compilation mode.
Columns {\tt LP-Dom} show that applying only the LP-dominance rule does not necessarily (highly) accelerate the DP algorithm. It is sometimes even more expensive to apply this version. That can be explained by the fact that the DP algorithm has to be slightly adapted to deal with classes in the problem with potentially 1, 2 or 3 items. Thus the structure of the DKP is disrupted and additional tests are needed to apply the recurrence. However, a very interesting result come from column {\tt opt-lb}: 
this (heuristic) fixation discard only one optimal solution for instance {\tt sdkp1\_6} in the first set and only for four instances in set 2. 
Using this heuristic fixation can be a very interesting alternative since even when the optimal solutions are discarded the difference between the optimal values is almost negligible (i.e., between 1 and 5 units). Please note that we do not provide the results obtained when using only the LP-dominance rule and looking for an optimal solution with the DP algorithm since the gain is too limited in that case (as in the case without solution). Thus these results do not provide additional information.

\begin{table}[htb]
\caption{Results obtained when combining DP algorithm and reduction.} \label{tab::dp}
\begin{center}
\begin{scriptsize}
\begin{tabular}{|I{2.2cm}|I{1cm}|I{1cm}|I{1cm}|I{1.3cm}|I{1cm}|I{1.1cm}|I{1.3cm}|I{1cm}|}
%
\cline{2-9}
\multicolumn{1}{c|}{} & \multicolumn{5}{c|}{\textbf{DP without solution}} & \multicolumn{3}{c|}{\textbf{DP with solution}} \\
\cline{2-9}
\multicolumn{1}{c|}{} & \textbf{Init.} & \multicolumn{2}{c|}{\textbf{LP-Dom}} & \textbf{Red. +} & \multirow{2}{*}{\textbf{Red.}} & \textbf{Init.} & \textbf{Red. +} & \multirow{2}{*}{\textbf{Red.}} \\
\cline{3-4}
\multicolumn{1}{c|}{} &  \textbf{Pb.}     & \textbf{CPU} & \textbf{opt-lb} & \textbf{LP-Dom}      &       &  \textbf{Pb.}     &  \textbf{LP-Dom}     & \\
\hline
Unc.  & 12.81 & 13.42 & 0     & 0.61  & 0.83  & 16.74 & 0.54  & 0.80 \\
\hline
Weak. & 4.88  & 3.82  & 0     & 0.14  & 0.20  & 10.24 & 0.16  & 0.22 \\
\hline
Strong. & 5.30  & 4.14  & 1     & 0.17  & 0.21  & 9.93  & 0.18  & 0.24 \\
\hline
Inv.  & 6.75  & 5.83  & 0     & 0.17  & 0.23  & 12.38 & 0.17  & 0.24 \\
\hline
Overall Group 1 & 7.44  & 6.80  & 1     & 0.27  & 0.37  & 12.21 & 0.26  & 0.38 \\
\hline
Unc.  & 1.18  & 1.05  & 2     & 0.16  & 0.17  & 1.19  & 0.16  & 0.19 \\
\hline
Weak. & 1.15  & 0.84  & 1     & 0.10  & 0.12  & 1.13  & 0.11  & 0.12 \\
\hline
Strong. & 1.46  & 0.96  & 0     & 0.12  & 0.13  & 1.08  & 0.15  & 0.16 \\
\hline
Inv.  & 1.26  & 1.03  & 1     & 0.07  & 0.10  & 1.16  & 0.09  & 0.09 \\
\hline
Overall Group 2 & 1.26  & 0.97  & 4     & 0.11  & 0.13  & 1.14  & 0.13  & 0.14 \\
\hline
Overall 80 inst. & 4.35  & 3.89  & 5     & 0.19  & 0.25  & 6.60  & 0.20  & 0.26 \\
\hline
\end{tabular}%
\end{scriptsize}
\end{center}
\end{table}

Columns {\tt Red.+LP-Dom} provide the running times needed by the DP algorithm when applying the class reduction technique and the LP-dominance fixation. In that case the value returned by the DP algorithm is necessarily the same as in the case of {\tt LP-Dom}. However here the class reduction process allows a very impressive speed-up of the DP algorithm since practically all the instances can be solved in less than one second with the backtracking phase to collect a solution. Thus this technique correspond to a very fast and effective heuristic to solve these two data sets of instances. 
Finally, columns {\tt Red.} give the running time needed by the DP algorithm to solve the reduced problem obtained when applying only the class reduction process. In that case we know that the value returned by the algorithm is the optimal value of the original DKP instance since this reduction is valid for the DKP. Here again values reported in these columns are quite impressive and demonstrate that solving exactly these two sets of instances can be done efficiently using this reduction rule. Detailed results showed that the running time with the backtracking phase is always less than 2 second and almost always less than 1 second. 

\section{Conclusion}
\label{sec::concl}

This paper deals with the discounted 0–1 knapsack problem (DKP), an extension of the knapsack problem, where items are grouped by three,  and at most one item from a group can be included in a solution. To tackle the problem we proposed two fixation techniques. The first one can be used in a heuristic way by discarding some variables of the problem according to a LP-dominance rule. This process can discard an optimal solution of the problem but it can be used to provide a near optimal feasible solution. The second one removes some groups from the problem, leading to a smaller problem to solve, without discarding any optimal solution of the original problem. Computational experiments were conducted to show that both fixation techniques can be used in a preprocessing phase before applying dynamic programming to solve the instances proposed in the literature in at most two seconds.
Future works will address the study of the structure of the instances of the DKP. Indeed experiments showed that the average gap between the linear programming relaxation and the optimal value of the problem is always very tight for the available instances. It could be interesting to consider other exact approaches based on branch-and-bound, for instance, and to look if some properties can be highlighted.


%

\begin{appendices}

\section{Detailed results: LP-relaxation and fixation}
\label{app:lp}

In all the appendixes instances in set 1 (resp. 2) have in their name {\tt "1\_"} (resp. {\tt "2\_")} before the number between 1 and 10. This prefix is used to denote the membership of instances.


We first provide in Table \ref{tab:detailed_lp_set1} and \ref{tab:detailed_lp_set2} the detailed results when solving the LP-relaxation of the instances and applying our fixation techniques. The columns meaning is the same as in Table \ref{tab:lp}: column {\tt LP-value} provides the LP-relaxation value whereas {\tt $\lfloor\text{LP}\rfloor$-opt} gives the difference between the value obtained by rounding down  LP value and the optimal value of the problem. Column {\tt opt-LPs} corresponds to the lower bound value obtained when extending the rounded solution of the LP-relaxation in a greedy way. Column {\tt LP-Dom} provides the percentage of variables that are LP-dominated. Column {\tt opt-lb} gives the difference between the optimal value and the final lower bound returned by the procedure. Columns {\tt Red.} and {\tt Red.+LP-Dom} refers to the percentage of classes that can be definitively eliminated from the problem and the total percentage of variables that can be set to 1 or 0 if we combine the two methods, respectively. Finally, columns {\tt CPU} gives the CPU time in seconds needed by the corresponding procedure.
\begingroup
\begin{center}
\begin{scriptsize}
\begin{longtable}{|C{1.2cm}|C{1.2cm}|C{1.5cm}|C{0.8cm}|C{0.8cm}|C{0.8cm}|C{0.8cm}|C{0.8cm}|C{0.8cm}|C{1.2cm}|C{0.8cm}|}
\caption{Results for the LP-relaxation and reduction process over set 1.} \label{tab:detailed_lp_set1} \\
\cline{3-11}
\multicolumn{2}{c}{} & \multicolumn{5}{|c|}{\textbf{LP solving results}} & \multicolumn{4}{c|}{\textbf{Reduction phase results}} \\
\hline
\textbf{Instance} & \textbf{opt} & \textbf{LP value} & \textbf{$\lfloor\text{LP}\rfloor$-opt} & \textbf{opt-LPg} & \textbf{LP-Dom} & \textbf{CPU} & \textbf{opt-lb} & \textbf{Red.} & \textbf{Red. + LP-Dom} & \textbf{CPU} \\
\hline
\endfirsthead
\multicolumn{11}{c}%
{{\tablename\ \thetable{} -- \textsl{continued from previous page}}} \\
\cline{3-11}
\multicolumn{2}{c}{} & \multicolumn{5}{|c|}{\textbf{LP solving results}} & \multicolumn{4}{c|}{\textbf{Reduction phase results}} \\
\hline
\textbf{Instance} & \textbf{opt} & \textbf{LP value} & \textbf{$\lfloor\text{LP}\rfloor$-opt} & \textbf{opt-LPg} & \textbf{LP-Dom} & \textbf{CPU} & \textbf{opt-lb} & \textbf{Red.} & \textbf{Red. + LP-Dom} & \textbf{CPU} \\
\hline
\endhead
\hline \multicolumn{11}{|r|}{{\textsl{continued on next page}}} \\ \hline
\endfoot
\endlastfoot
udkp1\_1 & 289761 & 289887.28 & 126   & 39    & 49.33 & 0.002 & 39    & 68.00 & 82.00 & 0.001 \\
\hline
udkp1\_2 & 510131 & 510184.73 & 53    & 145   & 47.67 & 0.002 & 103   & 49.50 & 70.50 & 0.003 \\
\hline
udkp1\_3 & 817713 & 817768.74 & 55    & 8     & 46.33 & 0.004 & 0     & 61.00 & 77.11 & 0.006 \\
\hline
udkp1\_4 & 1122074 & 1122122.50 & 48    & 177   & 46.50 & 0.003 & 45    & 64.25 & 79.17 & 0.014 \\
\hline
udkp1\_5 & 1233057 & 1233094.43 & 37    & 227   & 46.27 & 0.004 & 29    & 54.00 & 73.33 & 0.025 \\
\hline
udkp1\_6 & 1399458 & 1399484.24 & 26    & 273   & 46.56 & 0.003 & 20    & 50.33 & 71.67 & 0.036 \\
\hline
udkp1\_7 & 1826261 & 1826275.41 & 14    & 84    & 46.19 & 0.007 & 26    & 58.57 & 75.57 & 0.049 \\
\hline
udkp1\_8 & 1920409 & 1920432.39 & 23    & 196   & 46.96 & 0.005 & 31    & 51.25 & 72.00 & 0.057 \\
\hline
udkp1\_9 & 2458318 & 2458338.59 & 20    & 42    & 46.22 & 0.009 & 41    & 60.78 & 77.22 & 0.066 \\
\hline
udkp1\_10 & 2886506 & 2886532.21 & 26    & 213   & 46.93 & 0.005 & 81    & 69.50 & 82.13 & 0.134 \\
\hline
wdkp1\_1 & 310805 & 310906.38 & 101   & 780   & 65.33 & 0.001 & 0     & 60.00 & 85.33 & 0.001 \\
\hline
wdkp1\_2 & 504177 & 504217.25 & 40    & 346   & 64.83 & 0.002 & 126   & 45.50 & 80.00 & 0.003 \\
\hline
wdkp1\_3 & 840609 & 840641.51 & 32    & 580   & 63.56 & 0.002 & 12    & 60.00 & 83.56 & 0.006 \\
\hline
wdkp1\_4 & 1041019 & 1041063.46 & 44    & 9     & 64.58 & 0.006 & 9     & 52.75 & 82.17 & 0.012 \\
\hline
wdkp1\_5 & 1606341 & 1606380.50 & 39    & 9     & 64.20 & 0.002 & 9     & 69.60 & 87.40 & 0.044 \\
\hline
wdkp1\_6 & 1875732 & 1875753.93 & 21    & 163   & 65.00 & 0.007 & 13    & 68.00 & 87.67 & 0.055 \\
\hline
wdkp1\_7 & 1726671 & 1726687.66 & 16    & 294   & 65.10 & 0.004 & 35    & 50.71 & 82.00 & 0.05 \\
\hline
wdkp1\_8 & 2589429 & 2589446.22 & 17    & 226   & 64.63 & 0.012 & 35    & 68.88 & 87.58 & 0.123 \\
\hline
wdkp1\_9 & 2551957 & 2551984.69 & 27    & 60    & 64.41 & 0.003 & 10    & 60.78 & 84.67 & 0.105 \\
\hline
wdkp1\_10 & 2718419 & 2718434.37 & 15    & 340   & 64.53 & 0.007 & 31    & 55.80 & 83.13 & 0.108 \\
\hline
sdkp1\_1 & 352019 & 352051.56 & 32    & 1198  & 60.67 & 0.001 & 91    & 65.00 & 82.33 & 0.002 \\
\hline
sdkp1\_2 & 545255 & 545281.08 & 26    & 108   & 60.67 & 0.003 & 74    & 52.50 & 78.17 & 0.004 \\
\hline
sdkp1\_3 & 986019 & 986054.24 & 35    & 638   & 60.00 & 0.009 & 64    & 58.33 & 79.67 & 0.007 \\
\hline
sdkp1\_4 & 1247191 & 1247212.77 & 21    & 21    & 60.17 & 0.003 & 21    & 67.25 & 83.75 & 0.016 \\
\hline
sdkp1\_5 & 1759075 & 1759101.72 & 26    & 72    & 60.60 & 0.003 & 28    & 70.80 & 84.93 & 0.027 \\
\hline
sdkp1\_6 & 1795393 & 1795414.18 & 21    & 441   & 60.83 & 0.012 & 24    & 62.50 & 81.94 & 0.068 \\
\hline
sdkp1\_7 & 2264218 & 2264238.98 & 20    & 329   & 60.52 & 0.004 & 6     & 66.86 & 83.33 & 0.103 \\
\hline
sdkp1\_8 & 2236703 & 2236718.08 & 15    & 553   & 60.13 & 0.004 & 18    & 58.38 & 79.96 & 0.095 \\
\hline
sdkp1\_9 & 3034816 & 3034828.50 & 12    & 381   & 60.22 & 0.005 & 25    & 70.00 & 84.93 & 0.182 \\
\hline
sdkp1\_10 & 2916217 & 2916224.34 & 7     & 356   & 60.47 & 0.005 & 37    & 58.10 & 80.13 & 0.206 \\
\hline
idkp1\_1 & 277642 & 277729.12 & 87    & 458   & 63.33 & 0.002 & 0     & 51.00 & 80.33 & 0.001 \\
\hline
idkp1\_2 & 541724 & 541765.81 & 41    & 962   & 64.50 & 0.002 & 80    & 52.50 & 82.00 & 0.003 \\
\hline
idkp1\_3 & 1016524 & 1016575.83 & 51    & 597   & 65.22 & 0.002 & 40    & 65.67 & 87.11 & 0.008 \\
\hline
idkp1\_4 & 1220338 & 1220382.63 & 44    & 69    & 64.75 & 0.003 & 19    & 59.50 & 84.58 & 0.019 \\
\hline
idkp1\_5 & 1342480 & 1342525.58 & 45    & 305   & 65.47 & 0.004 & 37    & 51.80 & 82.73 & 0.026 \\
\hline
idkp1\_6 & 1922488 & 1922511.13 & 23    & 603   & 65.06 & 0.005 & 41    & 67.67 & 87.61 & 0.061 \\
\hline
idkp1\_7 & 2190780 & 2190809.56 & 29    & 178   & 64.95 & 0.004 & 33    & 66.14 & 87.00 & 0.061 \\
\hline
idkp1\_8 & 2719899 & 2719933.58 & 34    & 306   & 65.00 & 0.003 & 14    & 71.25 & 88.75 & 0.066 \\
\hline
idkp1\_9 & 2377631 & 2377643.96 & 12    & 303   & 64.74 & 0.005 & 50    & 50.44 & 81.56 & 0.126 \\
\hline
idkp1\_10 & 3123425 & 3123455.87 & 30    & 543   & 64.80 & 0.005 & 51    & 63.30 & 85.90 & 0.096 \\
\hline
\end{longtable}
\end{scriptsize}
\end{center}
\endgroup

\begingroup
\begin{center}
\begin{scriptsize}
\begin{longtable}{|C{1.3cm}|C{1.3cm}|C{1.5cm}|C{0.8cm}|C{0.8cm}|C{0.8cm}|C{0.8cm}|C{0.8cm}|C{0.8cm}|C{1.2cm}|C{0.8cm}|}
\caption{Results for the LP-relaxation and reduction process over set 2.} \label{tab:detailed_lp_set2} \\
\cline{3-11}
\multicolumn{2}{c}{} & \multicolumn{5}{|c|}{\textbf{LP solving results}} & \multicolumn{4}{c|}{\textbf{Reduction phase results}} \\
\hline
\textbf{Instance} & \textbf{opt} & \textbf{LP value} & \textbf{$\lfloor\text{LP}\rfloor$-opt} & \textbf{opt-LPg} & \textbf{LP-Dom} & \textbf{CPU} & \textbf{opt-lb} & \textbf{Red.} & \textbf{Red. + LP-Dom} & \textbf{CPU} \\
\hline
\endfirsthead
\multicolumn{11}{c}%
{{\tablename\ \thetable{} -- \textsl{continued from previous page}}} \\
\cline{3-11}
\multicolumn{2}{c}{} & \multicolumn{5}{|c|}{\textbf{LP solving results}} & \multicolumn{4}{c|}{\textbf{Reduction phase results}} \\
\hline
\textbf{Instance} & \textbf{opt} & \textbf{LP value} & \textbf{$\lfloor\text{LP}\rfloor$-opt} & \textbf{opt-LPg} & \textbf{LP-Dom} & \textbf{CPU} & \textbf{opt-lb} & \textbf{Red.} & \textbf{Red. + LP-Dom} & \textbf{CPU} \\
\hline
\endhead
\hline \multicolumn{11}{|r|}{{\textsl{continued on next page}}} \\ \hline
\endfoot
\endlastfoot
udkp2\_1 & 85740 & 85757.52 & 17    & 146   & 46.00 & 0.002 & 19    & 60.00 & 78.00 & 0.002 \\
\hline
udkp2\_2 & 163744 & 163775.43 & 31    & 13    & 44.50 & 0.005 & 13    & 47.50 & 68.17 & 0.003 \\
\hline
udkp2\_3 & 269393 & 269401.81 & 8     & 5     & 44.56 & 0.005 & 5     & 68.00 & 80.44 & 0.011 \\
\hline
udkp2\_4 & 347599 & 347605.76 & 6     & 11    & 46.25 & 0.004 & 0     & 61.00 & 76.92 & 0.015 \\
\hline
udkp2\_5 & 442644 & 442656.61 & 12    & 105   & 43.73 & 0.007 & 11    & 62.00 & 77.33 & 0.053 \\
\hline
udkp2\_6 & 536578 & 536585.58 & 7     & 79    & 44.94 & 0.005 & 19    & 64.83 & 78.83 & 0.077 \\
\hline
udkp2\_7 & 635860 & 635866.05 & 6     & 30    & 47.14 & 0.01  & 11    & 67.29 & 80.43 & 0.105 \\
\hline
udkp2\_8 & 650206 & 650213.17 & 7     & 31    & 46.00 & 0.005 & 12    & 58.25 & 75.08 & 0.085 \\
\hline
udkp2\_9 & 718532 & 718537.68 & 5     & 13    & 45.74 & 0.006 & 4     & 58.56 & 75.85 & 0.109 \\
\hline
udkp2\_10 & 779460 & 779464.19 & 4     & 12    & 45.27 & 0.006 & 5     & 50.80 & 70.80 & 0.121 \\
\hline
wdkp2\_1 & 83098 & 83113.56 & 15    & 15    & 58.00 & 0.003 & 15    & 57.00 & 78.00 & 0.001 \\
\hline
wdkp2\_2 & 138215 & 138222.05 & 7     & 2     & 57.83 & 0.002 & 2     & 55.00 & 77.67 & 0.004 \\
\hline
wdkp2\_3 & 256616 & 256621.13 & 5     & 126   & 57.33 & 0.005 & 6     & 66.00 & 82.11 & 0.014 \\
\hline
wdkp2\_4 & 315657 & 315663.03 & 6     & 118   & 58.67 & 0.005 & 10    & 61.50 & 81.25 & 0.016 \\
\hline
wdkp2\_5 & 428490 & 428492.83 & 2     & 77    & 58.73 & 0.007 & 9     & 65.20 & 82.00 & 0.04 \\
\hline
wdkp2\_6 & 466050 & 466054.45 & 4     & 0     & 58.89 & 0.004 & 0     & 63.33 & 81.61 & 0.049 \\
\hline
wdkp2\_7 & 547683 & 547686.89 & 3     & 123   & 58.33 & 0.023 & 5     & 63.71 & 81.90 & 0.075 \\
\hline
wdkp2\_8 & 576959 & 576961.99 & 2     & 40    & 59.08 & 0.004 & 7     & 58.88 & 80.54 & 0.086 \\
\hline
wdkp2\_9 & 650660 & 650664.21 & 4     & 21    & 58.56 & 0.005 & 11    & 56.89 & 78.81 & 0.106 \\
\hline
wdkp2\_10 & 678967 & 678971.13 & 4     & 8     & 58.77 & 0.014 & 5     & 54.50 & 77.93 & 0.139 \\
\hline
sdkp2\_1 & 94459 & 94480.26 & 21    & 464   & 54.00 & 0.002 & 23    & 66.00 & 80.33 & 0.001 \\
\hline
sdkp2\_2 & 160805 & 160814.92 & 9     & 141   & 53.50 & 0.003 & 4     & 61.50 & 78.17 & 0.004 \\
\hline
sdkp2\_3 & 238248 & 238254.23 & 6     & 236   & 52.89 & 0.002 & 27    & 58.00 & 75.44 & 0.01 \\
\hline
sdkp2\_4 & 340027 & 340034.57 & 7     & 274   & 53.42 & 0.002 & 8     & 62.00 & 78.00 & 0.019 \\
\hline
sdkp2\_5 & 463033 & 463034.99 & 1     & 372   & 52.53 & 0.006 & 28    & 65.40 & 79.27 & 0.071 \\
\hline
sdkp2\_6 & 466097 & 466100.89 & 3     & 65    & 50.89 & 0.003 & 22    & 57.00 & 73.89 & 0.054 \\
\hline
sdkp2\_7 & 620446 & 620451.26 & 5     & 142   & 52.24 & 0.004 & 20    & 65.14 & 78.90 & 0.094 \\
\hline
sdkp2\_8 & 670697 & 670702.05 & 5     & 352   & 53.08 & 0.005 & 15    & 61.13 & 77.88 & 0.121 \\
\hline
sdkp2\_9 & 739121 & 739126.94 & 5     & 284   & 52.63 & 0.005 & 11    & 60.89 & 77.33 & 0.136 \\
\hline
sdkp2\_10 & 765317 & 765319.32 & 2     & 338   & 51.83 & 0.005 & 2     & 59.00 & 76.43 & 0.181 \\
\hline
idkp2\_1 & 70106 & 70135.84 & 29    & 80    & 61.00 & 0.002 & 0     & 51.00 & 78.00 & 0.001 \\
\hline
idkp2\_2 & 118268 & 118275.19 & 7     & 36    & 58.50 & 0.002 & 0     & 48.50 & 74.67 & 0.002 \\
\hline
idkp2\_3 & 234804 & 234809.76 & 5     & 65    & 59.78 & 0.004 & 2     & 62.33 & 80.56 & 0.006 \\
\hline
idkp2\_4 & 282591 & 282599.51 & 8     & 26    & 60.25 & 0.003 & 0     & 58.50 & 79.75 & 0.012 \\
\hline
idkp2\_5 & 335584 & 335589.28 & 5     & 4     & 59.53 & 0.003 & 0     & 57.20 & 78.60 & 0.023 \\
\hline
idkp2\_6 & 452463 & 452467.40 & 4     & 18    & 60.06 & 0.003 & 0     & 64.00 & 81.39 & 0.032 \\
\hline
idkp2\_7 & 489149 & 489153.60 & 4     & 100   & 59.90 & 0.006 & 12    & 54.86 & 78.19 & 0.065 \\
\hline
idkp2\_8 & 533841 & 533850.27 & 9     & 25    & 59.71 & 0.004 & 24    & 48.63 & 75.92 & 0.056 \\
\hline
idkp2\_9 & 528144 & 528151.65 & 7     & 35    & 59.89 & 0.006 & 4     & 47.56 & 75.74 & 0.066 \\
\hline
idkp2\_10 & 581244 & 581250.49 & 6     & 19    & 59.67 & 0.05  & 6     & 46.40 & 75.13 & 0.09 \\
\hline
\end{longtable}
\end{scriptsize}
\end{center}
\endgroup

\newpage
\section{Detailed results: dynamic programming with fixation}
\label{app:dp}

In this section we provide in Tables \ref{tab::dp_det_set1} and \ref{tab::dp_det_set2} the detailed results obtained with the dynamic programming. The column meaning is the same as in Table \ref{tab::dp}: columns {\tt Init.Pb.} give the CPU times needed to solve the initial problem with the DP algorithm. Columns {\tt LP-Dom} report the running times of DP algorithm when it is applied on the reduced problem after applying the LP-dominance rule only and the difference between the optimal value and the value obtained by the DP algorithm. Column {\tt Red.+LP-Dom} provides the running time when we combine the class reduction phase with the LP-dominance. Columns {\tt Red.} contain the running times when the DP algorithm is applied on the reduced problem obtained when applying only the class reduction.

\begingroup
\begin{center}
\begin{scriptsize}
\begin{longtable}{|C{1.5cm}|C{1.3cm}|C{1.3cm}|C{1.3cm}|C{1.5cm}|C{1.3cm}|C{1.3cm}|C{1.3cm}|C{1.3cm}|}
\caption{Results obtained when combining DP algorithm and reduction over set 1.} \label{tab::dp_det_set1} \\
\cline{2-9}
\multicolumn{1}{c|}{} & \multicolumn{5}{c|}{\textbf{DP without solution}} & \multicolumn{3}{c|}{\textbf{DP with solution}} \\
\cline{2-9}
\multicolumn{1}{c|}{} & \textbf{Init.} & \multicolumn{2}{c|}{\textbf{LP-Dom}} & \textbf{Red. +} & \multirow{2}{*}{\textbf{Red.}} & \textbf{Init.} & \textbf{Red. +} & \multirow{2}{*}{\textbf{Red.}} \\
\cline{3-4}
\multicolumn{1}{c|}{} &  \textbf{Pb.}     & \textbf{CPU} & \textbf{opt-lb} & \textbf{LP-Dom}      &       &  \textbf{Pb.}     &  \textbf{LP-Dom}     & \\
\hline
\endfirsthead
\multicolumn{9}{c}%
{{\tablename\ \thetable{} -- \textsl{continued from previous page}}} \\
\cline{2-9}
\multicolumn{1}{c|}{} & \multicolumn{5}{c|}{\textbf{DP without solution}} & \multicolumn{3}{c|}{\textbf{DP with solution}} \\
\cline{2-9}
\multicolumn{1}{c|}{} & \textbf{Init.} & \multicolumn{2}{c|}{\textbf{LP-Dom}} & \textbf{Red. +} & \multirow{2}{*}{\textbf{Red.}} & \textbf{Init.} & \textbf{Red. +} & \multirow{2}{*}{\textbf{Red.}} \\
\cline{3-4}
\multicolumn{1}{c|}{} &  \textbf{Pb.}     & \textbf{CPU} & \textbf{opt-lb} & \textbf{LP-Dom}      &       &  \textbf{Pb.}     &  \textbf{LP-Dom}     & \\
\hline
\endhead
\hline \multicolumn{9}{|r|}{{\textsl{continued on next page}}} \\ \hline
\endfoot
\endlastfoot
udkp1\_1 & 0.337 & 0.191 & 0     & 0.031 & 0.033 & 0.223 & 0.038 & 0.046 \\
\hline
udkp1\_2 & 1.131 & 0.543 & 0     & 0.091 & 0.166 & 0.68  & 0.077 & 0.126 \\
\hline
udkp1\_3 & 2.971 & 2.263 & 0     & 0.155 & 0.229 & 1.934 & 0.143 & 0.21 \\
\hline
udkp1\_4 & 5.013 & 5.083 & 0     & 0.19  & 0.387 & 4.112 & 0.216 & 0.456 \\
\hline
udkp1\_5 & 7.116 & 5.333 & 0     & 0.327 & 1.036 & 8.602 & 0.402 & 0.694 \\
\hline
udkp1\_6 & 7.575 & 11.395 & 0     & 0.824 & 1.149 & 12.492 & 0.69  & 0.832 \\
\hline
udkp1\_7 & 14.087 & 16.968 & 0     & 0.988 & 1.29  & 27.385 & 0.748 & 1.074 \\
\hline
udkp1\_8 & 16.777 & 17.809 & 0     & 1     & 1.184 & 36.613 & 0.915 & 1.383 \\
\hline
udkp1\_9 & 28.273 & 30.949 & 0     & 1.079 & 1.559 & 58.612 & 1.224 & 1.794 \\
\hline
udkp1\_10 & 44.844 & 43.705 & 0     & 1.379 & 1.296 & O/M   & 0.974 & 1.38 \\
\hline
wdkp1\_1 & 0.216 & 0.151 & 0     & 0.008 & 0.011 & 0.162 & 0.015 & 0.016 \\
\hline
wdkp1\_2 & 0.317 & 0.282 & 0     & 0.034 & 0.044 & 0.395 & 0.04  & 0.065 \\
\hline
wdkp1\_3 & 1.035 & 0.711 & 0     & 0.02  & 0.027 & 1.102 & 0.037 & 0.045 \\
\hline
wdkp1\_4 & 1.353 & 0.945 & 0     & 0.072 & 0.093 & 1.837 & 0.085 & 0.115 \\
\hline
wdkp1\_5 & 3.275 & 1.824 & 0     & 0.096 & 0.096 & 6.36  & 0.091 & 0.109 \\
\hline
wdkp1\_6 & 4.669 & 3.017 & 0     & 0.118 & 0.144 & 9.223 & 0.116 & 0.147 \\
\hline
wdkp1\_7 & 4.223 & 3.225 & 0     & 0.196 & 0.309 & 8.849 & 0.276 & 0.34 \\
\hline
wdkp1\_8 & 9.762 & 7.142 & 0     & 0.304 & 0.327 & 20.524 & 0.261 & 0.381 \\
\hline
wdkp1\_9 & 10.387 & 9.473 & 0     & 0.236 & 0.317 & 22.891 & 0.32  & 0.424 \\
\hline
wdkp1\_10 & 13.603 & 11.428 & 0     & 0.356 & 0.59  & 31.011 & 0.34  & 0.566 \\
\hline
sdkp1\_1 & 0.236 & 0.096 & 0     & 0.013 & 0.016 & 0.196 & 0.018 & 0.019 \\
\hline
sdkp1\_2 & 0.426 & 0.25  & 0     & 0.035 & 0.035 & 0.524 & 0.052 & 0.062 \\
\hline
sdkp1\_3 & 1.412 & 0.708 & 0     & 0.088 & 0.103 & 1.288 & 0.126 & 0.164 \\
\hline
sdkp1\_4 & 2.272 & 1.253 & 0     & 0.066 & 0.061 & 2.297 & 0.07  & 0.086 \\
\hline
sdkp1\_5 & 2.682 & 2.62  & 0     & 0.092 & 0.113 & 4.19  & 0.166 & 0.179 \\
\hline
sdkp1\_6 & 3.255 & 2.301 & 5     & 0.15  & 0.191 & 4.073 & 0.167 & 0.243 \\
\hline
sdkp1\_7 & 5.18  & 3.411 & 0     & 0.187 & 0.213 & 12.435 & 0.213 & 0.276 \\
\hline
sdkp1\_8 & 5.257 & 4.342 & 0     & 0.264 & 0.267 & 13.943 & 0.337 & 0.41 \\
\hline
sdkp1\_9 & 15.165 & 11.527 & 0     & 0.355 & 0.436 & 29.376 & 0.312 & 0.385 \\
\hline
sdkp1\_10 & 17.159 & 14.883 & 0     & 0.479 & 0.708 & 31.014 & 0.384 & 0.559 \\
\hline
idkp1\_1 & 0.21  & 0.124 & 0     & 0.011 & 0.015 & 0.163 & 0.011 & 0.028 \\
\hline
idkp1\_2 & 0.653 & 0.282 & 0     & 0.018 & 0.035 & 0.535 & 0.026 & 0.036 \\
\hline
idkp1\_3 & 1.813 & 0.978 & 0     & 0.071 & 0.084 & 1.481 & 0.072 & 0.089 \\
\hline
idkp1\_4 & 2.486 & 1.389 & 0     & 0.089 & 0.126 & 2.255 & 0.103 & 0.158 \\
\hline
idkp1\_5 & 3.92  & 2.628 & 0     & 0.139 & 0.178 & 3.197 & 0.178 & 0.263 \\
\hline
idkp1\_6 & 5.365 & 3.441 & 0     & 0.157 & 0.192 & 6.265 & 0.136 & 0.241 \\
\hline
idkp1\_7 & 7.296 & 5.138 & 0     & 0.142 & 0.21  & 17.184 & 0.181 & 0.216 \\
\hline
idkp1\_8 & 11.674 & 10.735 & 0     & 0.172 & 0.262 & 26.169 & 0.214 & 0.296 \\
\hline
idkp1\_9 & 12.46 & 9.895 & 0     & 0.41  & 0.546 & 24.904 & 0.323 & 0.473 \\
\hline
idkp1\_10 & 21.581 & 23.642 & 0     & 0.527 & 0.685 & 41.627 & 0.462 & 0.645 \\
\hline
\end{longtable}
\end{scriptsize}
\end{center}
\endgroup

\begingroup
\begin{center}
\begin{scriptsize}
\begin{longtable}{|C{1.5cm}|C{1.3cm}|C{1.3cm}|C{1.3cm}|C{1.5cm}|C{1.3cm}|C{1.3cm}|C{1.3cm}|C{1.3cm}|}
\caption{Results obtained when combining DP algorithm and reduction over set 2.} \label{tab::dp_det_set2} \\
\cline{2-9}
\multicolumn{1}{c|}{} & \multicolumn{5}{c|}{\textbf{DP without solution}} & \multicolumn{3}{c|}{\textbf{DP with solution}} \\
\cline{2-9}
\multicolumn{1}{c|}{} & \textbf{Init.} & \multicolumn{2}{c|}{\textbf{LP-Dom}} & \textbf{Red. +} & \multirow{2}{*}{\textbf{Red.}} & \textbf{Init.} & \textbf{Red. +} & \multirow{2}{*}{\textbf{Red.}} \\
\cline{3-4}
\multicolumn{1}{c|}{} &  \textbf{Pb.}     & \textbf{CPU} & \textbf{opt-lb} & \textbf{LP-Dom}      &       &  \textbf{Pb.}     &  \textbf{LP-Dom}     & \\
\hline
\endfirsthead
\multicolumn{9}{c}%
{{\tablename\ \thetable{} -- \textsl{continued from previous page}}} \\
\cline{2-9}
\multicolumn{1}{c|}{} & \multicolumn{5}{c|}{\textbf{DP without solution}} & \multicolumn{3}{c|}{\textbf{DP with solution}} \\
\cline{2-9}
\multicolumn{1}{c|}{} & \textbf{Init.} & \multicolumn{2}{c|}{\textbf{LP-Dom}} & \textbf{Red. +} & \multirow{2}{*}{\textbf{Red.}} & \textbf{Init.} & \textbf{Red. +} & \multirow{2}{*}{\textbf{Red.}} \\
\cline{3-4}
\multicolumn{1}{c|}{} &  \textbf{Pb.}     & \textbf{CPU} & \textbf{opt-lb} & \textbf{LP-Dom}      &       &  \textbf{Pb.}     &  \textbf{LP-Dom}     & \\
\hline
\endhead
\hline \multicolumn{9}{|r|}{{\textsl{continued on next page}}} \\ \hline
\endfoot
\endlastfoot
udkp2\_1 & 0.032 & 0.035 & 0     & 0.006 & 0.006 & 0.055 & 0.011 & 0.015 \\
\hline
udkp2\_2 & 0.111 & 0.121 & 0     & 0.028 & 0.022 & 0.131 & 0.025 & 0.036 \\
\hline
udkp2\_3 & 0.361 & 0.398 & 0     & 0.032 & 0.031 & 0.304 & 0.025 & 0.036 \\
\hline
udkp2\_4 & 0.636 & 0.684 & 0     & 0.052 & 0.052 & 0.523 & 0.041 & 0.053 \\
\hline
udkp2\_5 & 1.046 & 1.217 & 2     & 0.128 & 0.126 & 0.883 & 0.094 & 0.153 \\
\hline
udkp2\_6 & 1.603 & 0.91  & 0     & 0.206 & 0.161 & 1.411 & 0.131 & 0.203 \\
\hline
udkp2\_7 & 1.639 & 1.404 & 0     & 0.214 & 0.196 & 1.702 & 0.153 & 0.239 \\
\hline
udkp2\_8 & 1.674 & 1.465 & 0     & 0.197 & 0.274 & 1.852 & 0.243 & 0.299 \\
\hline
udkp2\_9 & 1.598 & 1.772 & 0     & 0.258 & 0.239 & 2.283 & 0.256 & 0.3 \\
\hline
udkp2\_10 & 3.131 & 2.47  & 1     & 0.466 & 0.582 & 2.763 & 0.585 & 0.55 \\
\hline
wdkp2\_1 & 0.045 & 0.044 & 0     & 0.007 & 0.009 & 0.061 & 0.007 & 0.011 \\
\hline
wdkp2\_2 & 0.129 & 0.125 & 0     & 0.011 & 0.013 & 0.169 & 0.013 & 0.015 \\
\hline
wdkp2\_3 & 0.326 & 0.319 & 0     & 0.03  & 0.036 & 0.318 & 0.028 & 0.039 \\
\hline
wdkp2\_4 & 0.632 & 0.536 & 0     & 0.039 & 0.068 & 0.452 & 0.05  & 0.064 \\
\hline
wdkp2\_5 & 0.908 & 0.594 & 0     & 0.071 & 0.096 & 0.928 & 0.09  & 0.097 \\
\hline
wdkp2\_6 & 1.16  & 0.868 & 0     & 0.083 & 0.089 & 1.105 & 0.107 & 0.098 \\
\hline
wdkp2\_7 & 1.762 & 1.211 & 0     & 0.146 & 0.176 & 1.48  & 0.115 & 0.183 \\
\hline
wdkp2\_8 & 2.465 & 1.762 & 0     & 0.146 & 0.184 & 1.74  & 0.222 & 0.202 \\
\hline
wdkp2\_9 & 1.982 & 1.56  & 1     & 0.183 & 0.218 & 2.495 & 0.212 & 0.232 \\
\hline
wdkp2\_10 & 2.093 & 1.396 & 0     & 0.239 & 0.269 & 2.519 & 0.216 & 0.261 \\
\hline
sdkp2\_1 & 0.049 & 0.033 & 0     & 0.005 & 0.005 & 0.077 & 0.008 & 0.006 \\
\hline
sdkp2\_2 & 0.155 & 0.124 & 0     & 0.012 & 0.014 & 0.206 & 0.012 & 0.015 \\
\hline
sdkp2\_3 & 0.344 & 0.272 & 0     & 0.025 & 0.033 & 0.241 & 0.041 & 0.048 \\
\hline
sdkp2\_4 & 0.728 & 0.524 & 0     & 0.042 & 0.041 & 0.484 & 0.052 & 0.055 \\
\hline
sdkp2\_5 & 1.059 & 0.691 & 0     & 0.134 & 0.145 & 0.791 & 0.107 & 0.159 \\
\hline
sdkp2\_6 & 1.496 & 0.795 & 0     & 0.116 & 0.123 & 1.02  & 0.179 & 0.146 \\
\hline
sdkp2\_7 & 2.21  & 1.064 & 0     & 0.168 & 0.186 & 1.821 & 0.208 & 0.233 \\
\hline
sdkp2\_8 & 2.415 & 1.629 & 0     & 0.219 & 0.229 & 1.742 & 0.211 & 0.265 \\
\hline
sdkp2\_9 & 2.495 & 2.248 & 0     & 0.231 & 0.255 & 2.174 & 0.328 & 0.307 \\
\hline
sdkp2\_10 & 3.65  & 2.236 & 0     & 0.251 & 0.284 & 2.243 & 0.344 & 0.359 \\
\hline
idkp2\_1 & 0.048 & 0.037 & 0     & 0.006 & 0.005 & 0.061 & 0.005 & 0.008 \\
\hline
idkp2\_2 & 0.145 & 0.125 & 0     & 0.005 & 0.005 & 0.216 & 0.008 & 0.009 \\
\hline
idkp2\_3 & 0.534 & 0.378 & 1     & 0.02  & 0.021 & 0.393 & 0.024 & 0.027 \\
\hline
idkp2\_4 & 0.876 & 0.582 & 0     & 0.026 & 0.03  & 0.678 & 0.036 & 0.036 \\
\hline
idkp2\_5 & 0.904 & 0.721 & 0     & 0.035 & 0.038 & 0.92  & 0.04  & 0.037 \\
\hline
idkp2\_6 & 1.043 & 1.12  & 0     & 0.054 & 0.062 & 1.289 & 0.094 & 0.066 \\
\hline
idkp2\_7 & 2.171 & 1.245 & 0     & 0.103 & 0.148 & 1.403 & 0.115 & 0.133 \\
\hline
idkp2\_8 & 2.217 & 1.724 & 0     & 0.15  & 0.326 & 1.781 & 0.235 & 0.247 \\
\hline
idkp2\_9 & 2.379 & 1.981 & 0     & 0.133 & 0.149 & 2.085 & 0.183 & 0.146 \\
\hline
idkp2\_10 & 2.26  & 2.4   & 0     & 0.186 & 0.233 & 2.764 & 0.192 & 0.227 \\
\hline
\end{longtable}
\end{scriptsize}
\end{center}
\endgroup

\end{appendices}

\bibliographystyle{splncs04}

\bibliography{biblio}

\end{document}